 \documentclass{article}
\usepackage{amssymb}

\newtheorem{problem}{Problem}[section]
\newtheorem{theo}[problem]{Theorem}

\newtheorem{prop}[problem]{Proposition}
\newtheorem{cor}[problem]{Corollary}
\newtheorem{lema}[problem]{Lemma}

\begin{document}
\date{March 2003}
\title{{\Large Symmetric products of surfaces\\
and the cycle index}}

\author{{\Large Pavle Blagojevi\' c} \\
        {Mathematics Institute SANU, Belgrade} \\[2mm]
{\Large Vladimir Gruji\' c} \\ {Faculty of Mathematics,
Belgrade}\\[2mm]
  {\Large  Rade \v Zivaljevi\' c}
 \\ {Mathematics Institute
 SANU, Belgrade}}
\maketitle

\begin{abstract}

We study some of the combinatorial structures related to the
signature of $G$-symmetric products of (open) surfaces $SP^m_G(M)=
M^m/G$ where $G\subset S_m$. The attention is focused on the
question what information about a surface $M$ can be recovered
from a symmetric product $SP^n(M)$. The problem is motivated in
part by the study of locally Euclidean topological commutative
$(m+k,m)$-groups, \cite{TD2001}. Emphasizing a combinatorial point
of view we express the signature ${\rm Sign}(SP^m_G(M))$ in terms
of the cycle index $Z(G;\bar x)$ of $G$, a polynomial which
originally appeared in P{\' o}lya enumeration theory of graphs,
trees, chemical structures etc. The computations are used to show
that there exist punctured Riemann surfaces $M_{g,k}, M_{g',k'}$
such that the manifolds $SP^{m}(M_{g,k})$ and $SP^{m}(M_{g',k'})$
are often not homeomorphic, although they always have the same
homotopy type provided $2g+k = 2g'+k'$ and $k,k'\geq 1$.
\end{abstract}

\section{Introduction}
\label{intro}

The complex plane $\mathbb{C}$, the punctured plane
$\mathbb{C}^\ast = \mathbb{C}\setminus\{0\}$ and the elliptic
curves are classical examples of surfaces that support a group
structure. A natural generalization of the (commutative) group
structure is the structure of a (commutative) $(m+k,m)$-group.

Let $SP^n(X) := X^n/S_n$ be the symmetric product of $X$ and let
$SP^p(X)\times SP^q(X)\rightarrow SP^{p+q}(X), \, (a,b)\mapsto
c:=a\ast b$ be the operation induced by concatenation of strings
$a\in SP^p(X)$ and $b\in SP^q(X)$. A {\bf commutative}
$(m+k,m)$-groupoid is a pair $(X,\mu)$ where the
``multiplication'' $\mu$ is an arbitrary map $\mu :
SP^{m+k}(X)\rightarrow SP^m(X)$. The operation $\mu$ is
associative if for each $c\in SP^{m+2k}(X)$ and each
representation $c = a\ast b$, where $a\in SP^{m+k}(X)$ and $b\in
SP^m(X)$, the result $\mu(\mu(a)\ast b)$ is always the same, i.e.
independent from the particular choice of $a$ and $b$ in the
representation $c = a\ast b$. A commutative and associative
$(m+k,m)$-groupoid is a $(m+k,m)$-group if the equation $\mu(x\ast
a) = b$ has a solution $x\in SP^m(X)$ for each $a\in SP^k(X)$ and
$b\in SP^m(X)$. Note that $(2,1)$-groups are essentially the
groups in the usual sense of the word. If $X$ is a topological
space then $(X,\mu )$ is a topological $(m+k,m)$-group if it is a
$(m+k,m)$-group and the map $\mu : SP^{m+k}(X)\rightarrow SP^m(X)$
is continuous.

For the motivation and other information about commutative
$(m+k,m)$-groups the reader is referred to \cite{TD2001},
\cite{TD92}. Surprisingly enough, the only known surfaces that
support the structure of a $(m+k,m)$-group for $(m+k,m)\neq (2,1)$
are of the form $\mathbb{C}\setminus A$ where $A$ is a finite set.
Moreover it was proved in \cite{TD92}, see also Theorem~6.1 in
\cite{TD2001} that if $(M,\mu)$ is a locally Euclidean,
topological, commutative, $(m+k,m)$-group then $M$ must be an
orientable $2$-manifold. Moreover, a $2$-manifold that admits the
structure of a commutative $(m+k,m)$-group satisfies a strong
necessary condition that the symmetric power $SP^m(M):= M^m/S_m$
is of the form $\mathbb{R}^u\times (S^1)^v$. In particular the
signature ${\rm Sign}(SP^m(M))$ of $M$ must be zero.

There is a conjecture \cite{TD2001} that the only examples of
topological, commutative, $(m+k,m)$-groups are supported by
surfaces of the form $M = \mathbb{C}\setminus A$.
Corollary~\ref{sign-C} and Proposition~\ref{prop:wreath1} support
this conjecture, since they imply that all open surfaces of
sufficiently high genus have a non-zero signature. However this
conjecture serves also as a partial motivation for the following
general questions which may be of some independent interest.

\medskip\noindent
{\bf Questions:}

\begin{enumerate}
\item[(A)]
 To what extent is the topology of a surface $M$
determined by the topology of its symmetric product $SP^m(M)$ for
a given $m$?

\item[(B)] Are there examples of non-homeomorphic (open) surfaces
$M$ and $N$ such that the associated symmetric products $SP^m(M)$
and $SP^m(N)$ are homeomorphic?
\end{enumerate}

 In response to (A) we prove the following theorem
which says that homological invariants alone are not sufficient to
distinguish symmetric products of non-homeomorphic surfaces. This
puts some limitations on surfaces $M$ and $N$ in question (B),
however the question itself remains open and interesting already
in the case of general surfaces $M_{g,k}$.

\begin{theo}
\label{glavna} There exist open, orientable surfaces $M$ and $N$
such that the associated symmetric products $SP^{m}(M)$ and
$SP^{m}(N)$ are not homeomorphic although they have the same
homotopy type. More precisely, this is always true if $M =
M_{g,k}$ and $N = M_{g',k'}$ {\rm ($k,k'\geq 1$)} and
\begin{itemize}
\item $2g + k = 2g' + k'$, \item $g\neq g'$ and ${\rm
max}\{g,g'\}\geq m/2$
\end{itemize}
where $M_{g,k}:= M_g\setminus\{x_1,\ldots, x_k\}$ is the surface
of genus $g$ punctured at $k$ points.
\end{theo}

\medskip

A natural approach to questions ${\rm (A)}$ and ${\rm (B)} $ is to
determine what information about $M$ is hidden in $SP^m(M)$. If
$SP^m(M)$ is known, then the iterated symmetric product
$SP^k(SP^m(M))$ and its higher order analogs are also known. Since
$SP^k(SP^m(M))\cong SP_G(M) := M^{mk}/G$ where $G = S_k\wr S_m$ is
the {\em wreath} product of groups, it is natural to consider
general $G$-symmetric products $SP_G(M) = M^N/G$ where $G\subset
S_N$ is an arbitrary subgroup of $S_N$. When we want to emphasize
that $G$ is a subgroup of $S_N$ we write $SP_G(M) = SP^N_G(M)$.
$SP_G(M)$ is always a $\mathbb{Q}$-homology manifold and the
signature ${\rm Sign}(SP_G(M))$ is well defined. Our central
technical result is the following theorem.

\begin{theo}
\label{sign-E}
\begin{equation}
{\rm Sign}(SP^{m}_G(M_{g,k})) = Z(G;0,-2g,0, -2g,\ldots)
\end{equation}
where $G\subset S_m$ and $Z(G;x_1,x_2,\ldots, x_m)$ is the {\em
cycle index} of $G$.
\end{theo}

\begin{cor}
\label{sign-C}
\begin{equation}
\label{eqn:sign-3} {\rm Sign}(SP^{2n}(M_{g,k})) =  (-1)^n{g\choose
n}.
\end{equation}
\end{cor}

 Recall that the cycle index,
\cite{Aigner79}, \cite{Krishna}, \cite{GeSta}, \cite{Stanley}, is
a fundamental polynomial which was originally defined by
J.H.~Redfield and independently by G.~P{\' o}lya who recognized
its central role in his celebrated enumeration theory of groups,
graphs, trees, chemical compounds etc. The advantage of expressing
the result in terms of the cycle index lies in the fact that the
cycle index is a well studied object of enumerative combinatorics.
One of highlights is Theorem~\ref{Polya}, see \cite{Aigner79}
Section~V, which describes a procedure how the cycle index $Z(G;
\vec{x})$ of the wreath product $G = Q\wr H$ can be in a
transparent and elegant way expressed in terms of cycle indices
$Z(Q;\vec{x}), Z(H;\vec{x})$ of $Q$ and $H$ respectively. A good
illustration how these ideas can be applied is provided by the
following proposition which itself can be seen as another
corollary of Theorem~\ref{sign-E}.

\begin{prop}
\label{prop:wreath1} Suppose that $m$ is odd and $p$ an even
integer. Then
\begin{equation}
{\rm Sign}(SP_{S_p\wr S_m}(M_{g,k})) = Z(S_p\wr S_m; 0, -2g,
\ldots , 0, -2g) = (-1)^{p/2}{{\frac{1}{2}{2g\choose m}}\choose
\frac{1}{2}p} .
\end{equation}
\end{prop}

\begin{cor}
\label{sign-U} If $m$ is an odd integer and $k\geq 1$  then
\begin{equation}
\label{eqn:sign-4} {\rm Sign}(SP^2(SP^{m}(M_{g,k}))) =
-\frac{1}{2}{2g\choose m}.
\end{equation}
\end{cor}

For completeness we recall a remarkable formula of Don
Zagier\footnote{Zagier did not originally express the result in
terms of the cycle index $Z(G;\bar x)$. His result says that ${\rm
Sign}(\sigma_m,M^m)$ is either $\tau(M)$ or $\chi(M)$ depending on
the parity of $m$.} obtained by the Atiyah-Singer $G$-signature
theorem applied to $(\sigma_m,M^m)$, where $\sigma :
M^m\rightarrow M^m$ is the cyclic permutation.

\begin{theo}{\rm (Zagier \cite{Zagier1972})}
\label{Zagi}
\begin{equation}
\label{zagi-1} \begin{array}{rcl} {\rm Sign}(SP^m_G(M)) & = &
Z(G;\tau,\chi,\tau,\chi, \ldots \,)\\ & = & Z(G;\vec{x})\vert_{
x_{2i}=\chi,\, x_{2i-1}=\tau}
\end{array}
\end{equation}
where  $\tau = \tau(M)$ and $\chi = \chi(M)$ are respectively the
signature and the Euler characteristic of a compact, oriented,
even dimensional manifold without boundary.
\end{theo}

Again the use of the cycle index is convenient.  In the case $G =
S_k\wr S_m$ one obtains, along the lines of
Proposition~\ref{prop:wreath1}, the following result which for
$m=1$ reduces to the formula of Hirzebruch, \cite{Hirze1968},
\cite{Zagier1972}. By convention \cite{Stanley}, $[t^p]f(t)$ is
the coefficient of $t^p$ in the power series $f(t)$.

\begin{theo}
\label{Hirze}
\begin{equation}
\label{hirze-1} {\rm Sign}(SP_{S_k\wr S_m}(M)) =
[t^k](1-t^2)^{\frac{1}{2}(-1)^{m+1}{-\chi\choose
m}}\left(\frac{1+t}{1-t}\right)^{\frac{1}{2}{{\rm Sign}(SP^m(M))}}
\end{equation}
\end{theo}

\subsection{Signature as a function of both  $G$ and $g$}
\label{Gg-dependance}

Our proof of Theorem~\ref{sign-E} with minor modification yields a
proof of the corresponding well known statement for closed
surfaces $M_g$. This allows us to check our computations so we
find it convenient to formulate and prove these two results as
parts of a single statement.

\begin{theo}
\label{sign-A}
\begin{equation}
\label{sign-1}
\begin{array}{cccl}
{\rm (a)} & {\rm Sign}(SP^m_G(M_g)) &  = & Z(G;0,2-2g,0,2-2g,\ldots \,)\\
&& = & Z(G;\bar{x})\vert_{ x_{2i}=2-2g,\, x_{2i-1}=0}\\
{\rm (b)} & {\rm Sign}(SP^m_G(M_{g,k})) &  = & Z(G;0,-2g,0,-2g,\ldots \,)\\
&& = & Z(G;\bar{x})\vert_{ x_{2i}=-2g,\, x_{2i-1}=0}
\end{array}
\end{equation}
where $Z(G;x_1,\ldots ,x_m)$ is the cycle index of the permutation
group $G\subset S_m$.
\end{theo}

Before we commence the proof, let us recall some generalities
about the $g$-signature of $G$-manifolds or vector spaces with
$G$-invariant bilinear forms, \cite{HirzZagi}, \cite{Gord}.

Let $V$ be a vector space, $B : V\times V\rightarrow \mathbb{C}$ a
hermitian bilinear form on $V$ and $g : V\rightarrow V$ an
endomorphism which preserves the form $B$, $\, B(gx,gy) = B(x,y)$
for all $x,y \in V$. Then $V$ admits a $g$-invariant decomposition
$V \cong V^0\oplus V^+\oplus V^-$ such that $B$ is positive
definite on $V^+$, negative definite on $V^-$ and zero on $V^0$.
Then the $g$ signature of the triple $(V,B,g)$, or for short the
$g$-signature of $V$ is defined by
\[
\mbox{\rm Sign}(g,V) = \mbox{\rm Sign}(g,V,B) = \mbox{\rm
Trace}(g\vert V^+) - \mbox{\rm Trace}(g\vert V^-)
\]
A symmetric or skew-symmetric form $B : V\times V\rightarrow
\mathbb{R}$ defined on a real vector space $V$ can be extended to
a hermitian form $\widehat B$ on the complexified space $V\otimes
\mathbb{C}$ by the formula \cite{Gord},

\begin{equation}
\label{formula1} \widehat{B}(x\otimes\alpha,y\otimes\beta) =
\left\{
\begin{array}{rl}\alpha\bar\beta B(x,y), & B \mbox{ {\rm is symmetric}}\\
i\alpha\bar\beta B(x,y), & B \mbox{ {\rm is skew-symmetric}}
\end{array}\right.
\end{equation}
The associated $g$-signature is also denoted by $\mbox{\rm
Sign}(g,V)$. Finally, suppose that $M^{2n}$ is a smooth, oriented
manifold, with or without boundary, and let $g : M^{2n}\rightarrow
M^{2n}$ be an orientation preserving diffeomorphism of $M$. The
intersection form $B : H_n(M,\mathbb{Q})\times
H_n(M,\mathbb{Q})\rightarrow \mathbb{Q}$ is symmetric if $n$ is
even, or skew-symmetric if $n$ is odd, and the associated
endomorphism $\tilde{g} := H_n(g) : H_n(M)\rightarrow H_n(M)$
preserves both the intersection form $B$ and its complexification.
The $g$-signature of $(\tilde{g},V,B)$ is in this case denoted by
$\mbox{\rm Sign}(g,M)$. In particular, we observe that the usual
signature $\mbox{\rm Sign}(M)$ can be interpreted as the
$g$-signature $\mbox{\rm Sign}({\rm Id},M)$ of the identity map
${\rm Id} : M\rightarrow M$.

\medskip
The following well known result, \cite{Grot}, \cite{HirzZagi}, is
of fundamental importance. Note that even the case of
$0$-dimensional manifolds (finite sets) is interesting, when this
result reduces to an elementary lemma (Burnside lemma) which is a
corner-stone of P{\' o}lya enumeration theory.

\begin{prop}
\label{corner} Suppose that $M^{2n}$ is a smooth, oriented
manifold with a not necessarily free, orientation preserving
action of a finite group $G$ of diffeomorphism. More generally, it
is sufficient to assume that $M$ is a $\mathbb{Q}$-homology
manifold. Then $M/G$ is a $\mathbb{Q}$-homology manifold, ${\rm
Sign}(M/G)$ is well defined and the following formula holds,
{\rm\cite{Grot}, \cite{HirzZagi},}
\[
\mbox{\rm Sign}(M/G) = \frac{1}{\vert G\vert}\sum_{g\in
G}~\mbox{\rm Sign}(g,M).
\]
\end{prop}

The following proposition is used in the proof of
Theorem~\ref{sign-A}.

\begin{prop}
\label{cyclic} Let $V$ be a $(2n)$-dimensional complex vector
space and let $B : V\times V\rightarrow \mathbb{C}$ be a hermitian
form. Suppose $\omega : V\rightarrow V$ is an endomorphism
preserving the form $B$, such that $\omega^{2n}=1$, and for some
$v_0\in V$, the set
$\{v_0,\omega(v_0),\ldots,\omega^{2n-1}(v_0)\}$ is a basis of $V$.
Suppose that $B(v_0,\omega^j(v_0))\neq 0 \Longrightarrow j=n$.
Then,
\begin{equation}
\label{cyclic-1} \mbox{\rm Sign}(\omega, V) =
(\sum_{x^{2n}=1}~x^{n+1})\, {sign}\,B(v_0,\omega^n(v_0))
=\left\{\begin{array}{cl} 0, & n\neq 1\\ 2\,{sign}\,B(v_0,\omega
(v_0)), & n=1
\end{array}\right.
\end{equation}
\end{prop}

\noindent{\bf Proof:} Let
$x\in\{1,\epsilon,\epsilon^2,\ldots,\epsilon^{2n-1}\}$ be a root
of unity, $\epsilon = e^{2\pi i/2n}$. Let $z_x := v_0 +
x^{-1}\omega(v_0) + x^{-2}\omega^2(v_0)+\ldots
+x^{-(2n-1)}\omega^{2n-1}(v_0)$ be the eigenvector of $\omega$
corresponding to the eigenvalue $x$. If $x$ and $y$ a different
eigenvalues, $B(z_x,z_y) =
B(\omega(z_x),\omega(z_y))=x\bar{y}B(z_x,z_y)=0$. Otherwise, since
$B(v_0,\omega^j(v_0))$ can be nonzero only for $j=n$, we have
$B(z_x,z_x) = 2nx^nB(v_0,\omega^n(v_0))$. By a well known formula,
\[
{\rm Sign}(V,B) = \sum_{\lambda\in{\rm Spec}(\omega)}~\lambda\;
{\rm Sign}(V_{\lambda},B_{\lambda}),
\]
where $V_{\lambda}$ is the eigenspace of $\omega$ which
corresponds to the eigenvalue $\lambda$ and $B_{\lambda} = B\vert
V_{\lambda}$ is the restriction of the form $B$ on $V_{\lambda}$.
In our case,
\[
{\rm Sign}(V,B) = \sum_{x^{2n}=1} x\; {sign}\,B(z_x,z_x)=
(\sum_{x^{2n}=1}~x^{n+1}){sign}\,B(v_0,\omega^n(v_0))
\]
and the equation (\ref{cyclic-1}) follows.

\medskip

Our next step in the preparations for the proof of
Theorem~\ref{sign-A} is an explicit description of the
intersection pairing $B : H_{m}(M_g^m;\mathbb{Q})\times
H_{m}(M_g^m;\mathbb{Q}) \rightarrow \mathbb{Q}$. Note that both
sides of the equations (\ref{sign-1}) are zero if $m$ is an odd
number. So from here on, we focus our attention on the even case
and assume that $m=2n$.

Let us choose an orientation on $M_g$ and let $\mathbb{T}\in
H_2(M_g;\mathbb{Q})$ be the associated fundamental class. Let
$a_1, b_1, \ldots , a_g,b_g$ be a symplectic basis of $H_1(M_g)$
so that $a_i\cap a_i = b_j\cap b_j = 0, \, a_i\cap b_j= 0$ for
$i\neq j$, and $a_i\cap b_i = \mathbb{I}$, where $\mathbb{I}\in
H_0(M_g)$ is the generator. By K\" unneth formula,
\[
H_{2n}(M_g^{2n};\mathbb{Q}) \cong \bigotimes_{k_1+k_2\ldots
+k_{2n}=2n}~H_{k_i}(M_g;\mathbb{Q}).
\]
From here one deduces that a convenient basis of
$H_{2n}(M_g^{2n};\mathbb{Q})$ consists of all ``words'' $w =
w_1w_2\ldots w_{2n} = w_1\times w_2\times\ldots\times w_{2n}$
where $w_k\in \{a_i\}_{i=1}^g\cup
\{b_j\}_{j=1}^g\cup\{\mathbb{I},\mathbb{T}\}$. Note that for
dimensional reasons, the number of occurrences of the letter
$\mathbb{I}$ in the word $w$ is equal to the number of occurrences
of the letter $\mathbb{T}$.

\begin{lema}
\label{word} Let $w = w_1w_2\ldots w_{2n}$ and $w' =
w_1'w_2'\ldots w_{2n}'$ be two words, representing basic homology
classes in $H_{2n}(M_g^{2n})$. Then,
\begin{equation}
B(w,w') = \epsilon_{w,w'} \langle w_1 , w_1'\rangle \cdots \langle
w_{2n} , w_{2n}'\rangle
\end{equation}
where $\epsilon_{w,w'}$ is either $+1$ or $-1$, while
$B(\cdot,\cdot)$ and $\langle\cdot ,\cdot\rangle $ are the
intersection pairings on groups $H_{2n}(M_g^{2n})$ and $H_1(M_g)$
respectively.
\end{lema}

\noindent{\bf Proof:} Indeed, $w = w_1\times\ldots \times w_{2n} =
{\widehat{w}}_1\cap \ldots\cap\,\widehat{w}_{2n}$ where,
$\widehat{w}_i = \mathbb{T}\times\ldots\times
w_i\times\ldots\times\mathbb{T}\in H_{\ast}(M_g^{2n})$. Hence,
$B(w_1\times\ldots \times w_{2n}, w_1'\times\ldots \times w_{2n}')
= ({\widehat{w}'}_1\cap \ldots\cap\,\widehat{w}'_{2n})\cap
({\widehat{w}}_1\cap \ldots\cap\,\widehat{w}_{2n})$ and it is
sufficient to remember that
\[
a\cap b = (-1)^{{\rm cd}(a)\,{\rm cd}(b)}~b\cap a
\]
where ${\rm cd}(x)$ is the codimension of a class $x$.

\bigskip
Let us define an involution $\ast : H_{\ast}(M_g;\mathbb{Q})
\rightarrow H_{\ast}(M_g;\mathbb{Q})$ by the formula
\begin{equation}
\begin{array}{cccc}
a_i^{\ast} = b_i & b_j^{\ast} = a_j & \mathbb{I}^{\ast} =
\mathbb{T} & \mathbb{T}^{\ast} = \mathbb{I}
\end{array}
\end{equation}
i.e. the involution $\ast$ is up to sign, the Poincar\' e duality
map. For a given word $w = w_1w_2\ldots w_{2n}$ let  $w^\ast =
w_1^{\ast}w_2^{\ast}\ldots w_{2n}^{\ast}$. Note that the first
part of the following proposition is just a reformulation of
Lemma~\ref{word}, while the second part gives a precise formula
for the sign function $\epsilon_{w,w'}$.

\begin{prop}
\label{prop-word} $$ B(w,w')\neq 0 \Longrightarrow w' = w^\ast $$
$$ B(w,w^\ast) = (-1)^{\alpha(w)\choose 2}\langle w_1 ,
w_1^{\ast}\rangle \cdots \langle w_{2n} , w_{2n}^{\ast}\rangle =
(-1)^{{\alpha(w)\choose 2}+\beta(w)} $$ where $\beta(w)$ is the
number of occurrences of letters $b_1,\ldots ,b_g$ in $w$, while
$\alpha(w)$ is the number of occurrences of both $a_i$ and $b_j$
in $w$. If $\mathbb{I}$ and $\mathbb{T}$ do not appear in $w$
whatsoever, then $B(w,w^\ast) = +1$.
\end{prop}

Let $\pi\in S_{2n}$ be a permutation and let $\alpha(\pi) =
1^{\alpha_1}2^{\alpha_2}\cdots {(2n)}^{\alpha_{2n}}$ be the
associated partition of $[2n] =  \{1,2,\ldots ,2n\}$. In light of
the well known equality ${\rm Sign}(g\times h,M\times N) = {\rm
Sign}(g,M)\,{\rm Sign}(h,N)$ \cite{Gord} \cite{HirzZagi},

\begin{equation}
\label{jednacina} {\rm Sign}(\pi, (M_g)^{2n}) =
\prod_{k=1}^{2n}~\left\{{\rm Sign}(C_k,
(M_g)^k)\right\}^{\alpha_k}
\end{equation}
where $C_k : (M_g)^k\rightarrow (M_g)^k$ is a cyclic permutation
of coordinates, $$C_k(x_1,x_2,\ldots ,x_k) =
(x_2,\ldots,x_k,x_1).$$

\begin{prop}
\label{vazno}
\[
{\rm Sign}(C_k,(M_g)^k) = \left\{\begin{array}{rl} 2-2g, & k
\mbox{ {\rm even}}\\ 0\quad \; , & k \mbox{ {\rm odd}}
\end{array}\right.
\]
\end{prop}

\noindent{\bf Proof:}  Let $V = H_k(M_g)\otimes\mathbb{C}$ and let
$B : V\times V \rightarrow \mathbb{C}$ be the hermitian form
obtained by formula (\ref{formula1}) from the intersection form on
$(M_g)^k$. Let $\omega : V\rightarrow V$ be the map induced by
$C_k$. As before, $\omega : V\rightarrow V$ is a $B$-preserving
endomorphism. If $V$ admits a $B$-orthogonal decomposition of the
form $V = V_1\oplus V_2$ where $V_1$ and $V_2$ are
$\omega$-invariant subspaces, then
\begin{equation}
\label{decomposition} {\rm Sign}(\omega,V) = {\rm
Sign}(\omega_1,V_1) + {\rm Sign}(\omega_2,V_2)
\end{equation}
where $\omega_i := \omega\vert V_i$ is the restriction of $\omega$
on $V_i$. Given a word $w = w_1\ldots w_{k}$, let $V_w = {\rm
span}\{w,\omega(w),\ldots,\omega^{k-1}(w)\}$ be the minimal
$\omega$-invariant subspace of $V$ which contains vector $w$. Let
$W_w := V_w + V_{w^\ast}$. Then by Proposition~\ref{prop-word},
for any two words $w$ and $w'$, the associated spaces $W_w$ and
$W_{w'}$ are either identical or mutually orthogonal. By formula
(\ref{decomposition}), it suffices to find ${\rm
Sign}(\omega,W_w)$. If $V_w\neq V_{w^\ast}$ then ${\rm
Sign}(\omega,W_w) = 0$, so we assume that $W_w = V_w$. Note that
${\rm dim}(V_w) = p$, where $p = p(w):={\rm min}\{s\geq 1\mid
\omega^s(w) =\pm w\}$. If $k$ is odd, then $p$, being a divisor of
$k$, must be also an odd number. This implies that $w^\ast \neq
\pm\omega^j(w)$ for all $j$, which means that the form $B\vert
V_w$ is zero and ${\rm Sign}(\omega, V_w) = 0$. It immediately
follows that ${\rm Sign}(C_k,(M_g)^k) = 0$ if $k$ is an odd
number. If $k=2m$ is even, then a nonzero contribution from $V_w$
can be expected only if $p$ is an even number. In this case we can
apply Proposition~\ref{cyclic} and deduce that $p=2$. The list of
all basic words $u$ having the property $p(u)=2$ is $$ w =
\mathbb{I}\mathbb{T}\ldots\mathbb{I}\mathbb{T}, \qquad w_i:=
a_ib_i\ldots a_ib_i, \quad i=1,\ldots,g. $$ Since $\mathbb{I}$ and
$\mathbb{T}$ are classes of even degree, $\omega(w) = w^{\ast}$
and $B(w,\omega(w)) = +1$. Since $a_i, b_i$ are classes of degree
$1$, $B(w_i,\omega(w_i)) = (-1)^{2m-1}B(w_i,w_i^{\ast})$, and
knowing that by Proposition~\ref{prop-word},
$B(w_i,w_i^\ast)=(-1)^{\beta(w_i)+m(2m-1)}=+1$, we have
$B(w_i,\omega(w_i)) = -1$. Finally, by Proposition~\ref{cyclic},
\[
{\rm Sign}(C_{2m},(M_g)^{2m}) = {\rm Sign}(V_w) +
\sum_{i=1}^g~{\rm Sign}(V_{w_i}) = 2-2g.
\]

\begin{prop}
\label{vazno1} If $s\geq 1$ then
\[
{\rm Sign}(C_k,(M_{g,s})^k) = \left\{\begin{array}{rl} -2g, & k\,
\mbox{ {\rm even}}\\ 0, & k\, \mbox{ {\rm odd}} .
\end{array}\right.
\]
\end{prop}

\noindent{\bf Proof:} As in the proof of Proposition~\ref{vazno},
we need information about the intersection pairing
\[
B : H_{2n}(M_{g,k}^{\,2n};\mathbb{Q})\times
H_{2n}(M_{g,k}^{\,2n};\mathbb{Q}) \rightarrow \mathbb{Q}.
\]
The homology group $H_1(M_{g,s},\mathbb{Z})\cong
\mathbb{Z}^{2g+s-1}$ has a basis $e_1,e_2,\ldots, e_{2g+s-1}$
where $e_i = a_i, e_{g+j} = b_j$ for $i=1,\ldots,g$ and $e_j$ for
$j\geq 2g+1$ correspond to the holes $\alpha_1,\ldots
,\alpha_{s-1}$ in $M_g$. The group
$H_{2n}(M_{g,s}^{\,2n};\mathbb{Q})$ is generated by the classes
(words) of the form $w = w_1\times\cdots\times w_{2n} = w_1\cdots
w_{2n}$ where $w_i\in \{e_1,\ldots ,e_{2g+s-1}\}$.
Lemma~\ref{word} is still true with the simplification that both
$\mathbb{I}, \mathbb{T}$ and the classes
$\{e_j\}_{j=2g+1}^{2g+s-1}$, associated to the holes in $M_{g,s}$,
are excluded. The rest of the proof follows the argument of the
proof of Proposition~\ref{vazno}.

\medskip
\noindent{\bf Proof of Theorem~\ref{sign-A}:} As it was already
observed, Theorem~\ref{sign-A} is trivially correct if $m$ is an
odd number since in that case both sides of the equation
(\ref{sign-1}) are zero. Let us assume that $m$ is an even number,
$m=2n$. By equation (\ref{jednacina}) and Proposition~\ref{vazno}
\begin{equation}
\begin{array}{ccl}
{\rm Sign}(SP^{2n}_G(M_g)) & = & \frac{1}{(2n)!}\sum_{\pi\in
G}{\rm Sign}(\pi,(M_g)^{2n})\\ & = & \frac{1}{(2n)!}\sum_{\pi\in
G} \prod_{k=1}^{2n}~\left\{{\rm Sign}(C_k,
(M_g)^k)\right\}^{\alpha_k(\pi)}\\ & = & Z(G;0,2-2g,\ldots,0,2-2g)
\end{array}
\end{equation}
which establishes part ${\rm (a)}$ of Theorem~\ref{sign-A}. The
part ${\rm (b)}$ follows by the same argument from
Proposition~\ref{vazno1}.

\subsection{Corollaries of
Theorem~\ref{sign-E}}

We start with the following elementary lemma.

\begin{lema}
\label{elem}
\begin{equation}
\label{eqn:elem}
 Z(S_n;\alpha,\beta,\alpha,\beta,\ldots ) =
[t^n](1-t)^{-\frac{1}{2}(\alpha+\beta)}(1+t)^{\frac{1}{2}(\alpha-\beta)}
\end{equation}
\end{lema}

\medskip\noindent
{\bf Proof:} The result is easily deduced from the well known fact
\cite{Aigner79}, \cite{Krishna}, \cite{GeSta}, \cite{Stanley}
\begin{equation}
Z(S_n;x_1,\ldots, x_n) = [t^n]\exp (x_1t + \frac{x_2t^2}{2} +
\frac{x_3t^3}{3}+\ldots ).
\end{equation}

In particular for $\alpha=\beta$ one has
\begin{equation}
\label{jednakost} Z(S_n;\alpha,\alpha,\ldots ) =
[t^n](1-t)^{-\alpha} = (-1)^n{-\alpha\choose n}
\end{equation}
while for $\alpha = 0$ formula (\ref{eqn:elem}) reduces to

\begin{equation}
Z(S_n;0,\beta,0,\beta, \ldots ) = [t^n](1-t^2)^{-\beta/2}.
\end{equation}

\begin{theo}{\rm (G.~P{\' o}lya, \cite{Aigner79})}
\label{Polya} Let $G = S_k\wr H$ be the wreath product of $S_k$ by
a subgroup $H\subset S_m$ of the symmetric group $S_m$. Then the
cycle index $Z(G;\vec{x})$ of $G$ can be computed from the cycle
indices of $S_k$ and $H$ by the formula
\begin{equation}
Z(G;\vec{x}) = Z(S_k; Z(H;x_1,x_2,\ldots),
Z(H;x_2,x_4,\ldots),\ldots , Z(H;x_k,x_{2k},\ldots)).
\end{equation}
\end{theo}

\medskip\noindent
{\bf Proof of Proposition~\ref{prop:wreath1}:} By Theorems
\ref{sign-E} and \ref{Polya}
\[
\begin{array}{l}
{\rm Sign}(SP_{S_p\wr S_m}(M_{g,k})) \quad = \quad Z(S_p\wr S_m;
0, -2g,
0, -2g, \ldots )\\
= \quad Z(S_p; Z(S_m;0,-2g,0,-2g,\ldots), Z(S_m;
-2g,-2g,\ldots),\ldots ) \\
= \quad Z(S_p; 0, (-1)^m{2g\choose m}, 0, (-1)^m{2g\choose m},
\ldots ) \quad = \quad (-1)^{p/2}{{\frac{1}{2}{2g\choose
m}}\choose \frac{1}{2}p} .
\end{array}
\]

\subsection{Non-homeomorphic symmetric products}

\noindent {\bf Proof of Theorem~\ref{glavna}:} Since $M_{g,k}$ is,
up to homotopy, a wedge of $2g+k-1$ circles, the condition
$2g+k=2g'+k'$ implies $M_{g,k}\simeq M_{g',k'}$ and as a
consequence $SP^m(M_{g,k})\simeq SP^m(M_{g',k'})$. Suppose that
$m$ is even integer, $m=2n$. The open manifold $SP^{2n}(M_{g,k})$,
according to Corollary~\ref{sign-C}, has signature
$(-1)^n{g\choose n}$. The condition ${\rm max}\{g,g'\}\geq n$
guarantees that either ${g\choose n}$ or ${{g'}\choose n}$ is
nonzero. The sequence ${g\choose n}$, as a function of $g$, is
strictly monotone for $g\geq n$. Together with the condition
$g\neq g'$ this implies
\[
{\rm Sign}(SP^{2n}(M_{g,k})) \neq {\rm Sign}(SP^{2n}(M_{g',k'})),
\]
hence $SP^{2n}(M_{g,k})$ and $SP^{2n}(M_{g',k'})$ are not
homeomorphic. The case of an odd integer $m$ is treated similarly.
If contrary to the claim $SP^{m}(M_{g,k})$ and $SP^{m}(M_{g',k'})$
are homeomorphic, then by Corollary~\ref{sign-U} ${2g\choose m} =
{2g'\choose m}$. This again would contradict the conditions ${\rm
max}\{g,g'\}\geq m/2$ and $g\neq g'$.


\small \baselineskip3pt
\begin{thebibliography}{10}


\bibitem{Aigner79}
M.~Aigner.
\newblock {\em Combinatorial Theory,}
\newblock Springer-Verlag, Berlin 1979.





\bibitem{Bjo91}
A.~Bj{\"o}rner.
\newblock Topological methods.
\newblock In R.~Graham, M.~Gr{\"o}tschel, and L.~Lov{\'a}sz,
editors.
{\it   Handbook of Combinatorics}. North-Holland, Amsterdam, 1995.





\bibitem{GeSta}
I.M.~Gessel, R.P.~Stanley.
\newblock Algebraic enumeration.
\newblock In {\em Handbook of Combinatorics II}, R.L.~Grakham,
M.~Gr\" otschel, L.~Lov\' asz (eds.), North-Holland, Amsterdam
1995.




\bibitem{Gord}
C.M.~Gordon.
\newblock The $G$-signature theorem in dimension $4$,
\newblock in {\em A la Recherche de la Topologie Perdue},
L.~Guillou, A.~Marin (eds.), Progress in Math. 62, Birkh\" auser
1986.

\bibitem{GrHar}
P.~Griffiths, J.~Harris.
\newblock {\em Principles of Algebraic Geometry}.
\newblock John Wiley \& Sons, 1978.

\bibitem{Grot}
A.~Grothendieck.
\newblock Sur quelques points d'algebre homologique.
\newblock {\em T\^ ohoku Math. J.}, (2)\, 9 \, (1957), 119--221.


\bibitem{Hirze}
F.~Hirzebruch.
\newblock The signature theorem: reminiscences and recreation.
\newblock in {\em Prospects in mathematics,} Ann. of math. Studies
70, Princeton Univ. Press 1971, 3--31.

\bibitem{Hirze1968}
F.~Hirzebruch.
\newblock {\em Lectures on the Atiyah-Singer theorem and its
applications}.
\newblock (unpublished lecture notes), Berkeley 1968.

\bibitem{HirzZagi}
F.~Hirzebruch, D.~Zagier.
\newblock {\em The Atiyah -- Singer Theorem and Elementary Number
Theory}.
\newblock Publish or Perish, Inc., Boston, 1974.


\bibitem{Kallel_1998}
S.~Kallel.
\newblock Divisor spaces on punctured Riemann surfaces.
\newblock {\em Trans. Amer. Math. Soc.}, 350 (1998), 135--164.


\bibitem{Krishna}
V.~Krishnamurthy.
\newblock {\em Combinatorics, Theory and Applications,}
\newblock East--West Press 1985.

\bibitem{Macdonald_62}
I.G.~Macdonald.
\newblock The Poincar\' e polynomial of a symmetric product.
\newblock {\em Proc. Camb. Phil. Soc.}, 1962.

\bibitem{Macdonald_topo}
I.G.~Macdonald.
\newblock Symmetric product of an algebraic curve.
\newblock {\em Topology} 1 (1962), 319--343.


\bibitem{Stanley}
R.~Stanley
\newblock {\em Enumerative Combinatorics,} Vol. I and II.
\newblock Stud. Advanced Math. No. 49, Cambridge Univ. Press 2001.

\bibitem{TD92}
K.~Tren\v cevski, D.~Dimovski.
\newblock {\em Complex Commutative Vector Valued Groups.}
\newblock Macedonian Acad. Sci. and Arts, Skopje, 1992.

\bibitem{TD2001}
K.~Tren\v cevski, D.~Dimovski.
\newblock On the affine and projective commutative
$(m+k,m)$-groups.
\newblock {\em J. of Algebra}, 240 (2001), 338--365.


\bibitem{Zagier1972}
D.B.~Zagier.
\newblock {\em Equivariant Pontrjagin Classes and Applications to
Orbit Spaces.}
\newblock Lecture Notes in Mathematics 290, Springer-Verlag 1972.




\end{thebibliography}
\end{document}